\documentclass{dimacs-l}
\usepackage{amssymb,amsmath,latexsym}

\newtheorem{theorem}{Theorem}[section]
\newtheorem{corollary}[theorem]{Corollary}

\theoremstyle{definition}
\newtheorem{definition}[theorem]{Definition}

\numberwithin{equation}{section}

\def\N{\mathbb{N}}
\def\Q{\mathbb{Q}}
\def\R{\mathbb{R}}
\def\Z{\mathbb{Z}}

\def\K{{\mathcal K}}
\def\P{{\mathcal P}}

\def\c{{\bf c}}
\def\m{{\bf m}}
\def\x{{\bf x}}
\def\s{s}
\newcommand\const{\operatorname{const}} 

\begin{document}

\title{Dedekind sums: a combinatorial-geometric viewpoint} 

\author{Matthias Beck}
\address{Department of Mathematics, San Francisco State University, San Francisco, CA 94132, USA}
\email{beck@math.sfsu.org}

\author{Sinai Robins}
\address{Department of Mathematics, Temple University, Philadelphia, PA 19122, USA} 
\email{srobins@math.temple.edu} 
\thanks{The second author is supported by the NSA Young Investigator Grant MSPR-OOY-196.} 

\copyrightinfo{2004}
  {American Mathematical Society}

\subjclass[2000]{Primary 05A15, 11L03; Secondary 11P21, 52C07}
\date{April 24, 2000 and, in revised form December 6, 2001.} 

\keywords{Dedekind sums, rational polytopes, lattice points, partition function} 

\thanks{Appeared in: M.~B.~Nathanson (ed.), Unusual Applications of Number Theory, \emph{DIMACS: Series in Discrete Mathematics and Theoretical Computer Science} {\bf 64} (2004), 25--35.}

\begin{abstract}
The literature on Dedekind sums is vast.  In this expository
paper we show that there is
a common thread to many generalizations of Dedekind sums, namely through the 
study of lattice point enumeration of rational polytopes.   
In particular, there are some natural finite Fourier series which we call
{\bf Fourier-Dedekind sums}, and which form the building blocks of the number
of partitions of an integer from a finite set of positive integers.   This problem
also goes by the name of the `coin exchange problem'.
Dedekind sums have enjoyed a resurgence of interest recently, from such diverse
fields as 
topology, number theory, and combinatorial geometry.
The Fourier-Dedekind sums we study here include as special cases  generalized 
Dedekind sums studied by Berndt, Carlitz, Grosswald, Knuth, Rademacher, and
Zagier.   Our interest in these sums stems from 
the appearance of Dedekind's and Zagier's sums in lattice point count formulas 
for polytopes. 
Using some simple generating functions, we show that generalized Dedekind sums 
are natural 
ingredients for such formulas. As immediate `geometric' corollaries to our formulas, 
we obtain and generalize reciprocity laws of Dedekind, 
Zagier, and Gessel.   Finally, we prove a polynomial-time
complexity result for Zagier's higher-dimensional Dedekind sums. 
\end{abstract}

\maketitle


\section{Introduction} 
In recent years, Dedekind sums and their various siblings have enjoyed a new renaissance. 
Historically, they appeared in analytic number theory (Dedekind's $\eta$-function \cite{dedekind}), 
algebraic number theory (class number formulae \cite{meyer}), topology (signature defects of manifolds \cite{hirz}), 
combinatorial geometry (lattice point enumeration \cite{mordell}), and algorithmic complexity 
(pseudo random number generators \cite{knuth}).
In this expository paper, we define some broad generalizations of Dedekind sums, 
which are in fact 
finite Fourier series. We show that they appear naturally in the enumeration of 
lattice points in polytopes, 
and prove reciprocity laws for them. 

In combinatorial number theory, one is interested in partitions of an integer $n$
from a finite set. 
That is, one writes $n$ as a nonnegative integer linear combination of 
a given finite set of  positive
integers. We showed in \cite{bdr} that the number of such partitions of $n$ from 
a finite set 
is a quasipolynomial in $n$, whose coefficients are built up from the 
following generalization of Dedekind sums. 
\begin{definition} For $a_{ 0 } , \dots, a_{ d }, n \in \Z$, we define the {\bf Fourier-Dedekind sum} as
  \[ \sigma_{n} \left( a_{ 1 } , \dots, a_{ d } ; a_{0} \right) := \frac{1}{a_{0}} \sum_{ \lambda^{ a_{0} } = 1 } \frac{ \lambda^{ n } }{ \left( 1 - \lambda^{ a_{ 1 } } \right) \cdots \left( 1 - \lambda^{ a_{ d } } \right) } \ . \]
Here the sum is taken over all $a_{0}$'th roots of unity for which the summand is not singular.
\end{definition}
In \cite{gessel}, Gessel systematically studied sums of the form 
  \[ \sum_{ \lambda^{a} = 1 } R( \lambda ) \ , \] 
where $R$ is a rational function, and the sum is taken over all $a$'th 
roots of unity for which 
$R$ is not singular. He called them `generalized Dedekind sums', since his 
definition includes 
various generalizations of the Dedekind sum as special cases.
Hence we study Gessel's sums where the poles of $R$ are restricted to be roots of unity. 

In Section \ref{introded}, we give a brief history on those generalizations of 
the classical Dedekind 
sum (due to Rademacher \cite{rademacher}, and Zagier \cite{zagier}) which can be written as 
Fourier-Dedekind sums. Our interest in these sums stems from 
the appearance of Dedekind's and Zagier's sums in lattice point enumeration formulas for polytopes 
\cite{mordell,pommersheim,brion,diaz}. 
Using generating functions, we show in Section \ref{comb} that generalized Dedekind sums are natural ingredients 
for such formulas, which also apply to the theory of partition functions. 
In Section \ref{recsec} we obtain and generalize 
reciprocity laws of Dedekind \cite{dedekind}, Zagier \cite{zagier}, and Gessel \cite{gessel} 
as `geometric' corollaries to our formulas.   Finally, in Section \ref{last}, we prove
that Zagier's higher-dimensional Dedekind sums are in fact polynomial-time computable in fixed
dimension.   For Dedekind sums in 2 dimensions, this fact follows easily from their 
reciprocity law; but for higher dimensional Dedekind sums the polynomial-time
complexity does not seem to follow so easily, and we therefore invoke some recent work
of \cite{barv} and \cite{diaz}. 


\section{Classical Dedekind sums and generalizations}\label{introded} 
According to Riemann's will, it was his wish that Dedekind should get Riemann's unpublished 
notes and manuscripts \cite{grosswald}. Among these was a discussion of the 
important function 
  \[ \eta (z) = e^{ \frac{ \pi i z }{ 12 } } \prod_{ n \geq 1 } \left( 1 - e^{ 2 \pi i n z } \right) \ , \] 
which Dedekind took up and eventually published in Riemann's collected works \cite{dedekind}. 
\begin{definition} Let $ ((x)) $ be the sawtooth function defined by 
  \[ ((x)) := \left\{ \begin{array}{cl} \{ x \} - \frac{ 1 }{ 2 } & \mbox{ if } x \not\in \Z \\ 
                                        0                         & \mbox{ if } x \in \Z \ . \end{array} \right. \] 
Here $ \{ x \} = x - [x] $ denotes the fractional part of $x$. 
For two integers $a$ and $b$, we define the {\bf Dedekind sum} as 
  \[ \s (a,b) := \sum_{ k \text{ mod } b } \left( \left( \frac{ ka }{ b } \right) \right) \left( \left( \frac{ k }{ b } \right) \right) \ . \]
Here the sum is over a complete residue system modulo $b$. 
\end{definition} 
Through the study of the transformation properties of $ \eta $ under 
$ \mbox{SL}_{2} ( \Z ) $, Dede\-kind naturally arrived at $ \s (a,b) $. 
The classic introduction to the arithmetic properties of the Dedekind sum is \cite{grosswald}. 
The most important of these, already proved by Dedekind \cite{dedekind}, is the famous reciprocity law: 
\begin{theorem}[Dedekind]\label{dedrec} If $a$ and $b$ are relatively prime then 
  \[ \s (a,b) + \s (b,a) = - \frac{1}{4} + \frac{1}{12} \left( \frac{ a }{ b } + \frac{ 1 }{ ab } + \frac{ b }{ a } \right) \ . \] 
\end{theorem} 
This reciprocity law is easily seen to be equivalent to the transformation law of the $\eta$-function \cite{dedekind}. 
Due to the periodicity of $ (( x )) $, we can reduce $a$ modulo $b$ in the Dedekind sum:  
$ \s (a,b) = \s ( a \mbox{ mod } b , b ) $. Therefore, Theorem \ref{dedrec} allows us to 
compute $ \s (a,b) $ in polynomial time, similar in spirit to the Euclidean algorithm. 

The Dedekind sum $ \s (a,b) $ has various generalizations, two of which we introduce here. 
The first one is due to Rademacher \cite{rademacher}, who generalized sums introduced by 
Meyer \cite{meyer} and Dieter \cite{dieter}:  
\begin{definition} For $ a,b \in \Z $, $ x,y \in \R $, the {\bf Dedekind-Rademacher sum} is defined by 
  \[ \s (a,b;x,y) := \sum_{ k \text{ mod } b } \left( \left( \frac{ (k+y) a }{ b } + x \right) \right) \left( \left( \frac{ k+y }{ b } \right) \right) \ . \]
\end{definition} 
This sum posesses again a reciprocity law: 
\begin{theorem}[Rademacher]\label{radrec} If $a$ and $b$ are relatively prime and $x$ and $y$ are not both integers, then 
  \[ \s (a,b;x,y) + \s (b,a;y,x) = ((x)) ((y)) + \frac{ 1 }{ 2 } \left( \frac{ a }{ b } B_{ 2 } (y) + \frac{ 1 }{ ab } B_{ 2 } (ay+bx) + \frac{ b }{ a } B_{ 2 } (x) \right) . \] 
Here 
  \[ B_{ 2 } (x) := ( x - [x] )^{ 2 } - ( x - [x] ) + \frac{ 1 }{ 6 }   \]
is the periodized second Bernoulli polynomial. 
\end{theorem} 
If $x$ and $y$ are both integers, the Dedekind-Rademacher sum is simply the classical Dedekind sum, whose reciprocity law we already stated. 
As with the reciprocity law for the classical Dedekind sum, Theorem \ref{radrec} can be used to compute $ \s (a,b;x,y) $ in polynomial time.

The second generalization of the Dedekind sum we mention here is due to Zagier \cite{zagier}. 
From topological considerations, he arrived naturally at expressions
of the following kind:
\begin{definition}
Let $ a_{1} , \dots , a_{d} $ be integers relatively prime to $ a_{0} \in \N $. 
Define the 
{\bf higher-dimensional Dedekind sum} as 
  \[ \s ( a_{0} ; a_{1} , \dots , a_{d} ) := \frac{ (-1)^{d/2} }{ a_{0} } 
\sum_{ k=1 }^{ a_{0} - 1 } \cot \frac{ \pi k a_{1} }{ a_{0} } 
\cdots \cot \frac{ \pi k a_{d} }{ a_{0} } \ . \]
\end{definition} 
This sum vanishes if $d$ is odd. It is not hard to see that this indeed generalizes 
the classical Dedekind sum: the latter can be written in terms of cotangents 
\cite{grosswald}, which yields 
  \[ \s (a,b) = \frac{ 1 }{ 4b } \sum_{ k \text{ mod } b } 
\cot \frac{ \pi k a }{ b } \cot \frac{ \pi k }{ b } = - \frac{1}{4} \s ( b; a, 1 ) \ . \]
Again, there exists a reciprocity law for Zagier's sums: 
\begin{theorem}[Zagier]\label{zagierthm} If $ a_{0} , \dots , a_{d} \in \N$ 
are pairwise relatively prime then 
  \[ \sum_{ j=0 }^{ d } \s ( a_{j} ; a_{0} , \dots , \hat{ a_{j} } , \dots , a_{d} ) = 
\phi ( a_{0} , \dots , a_{d} ) \ . \] 
Here $ \phi $ is a rational function in $ a_{0} , \dots , a_{d} $, which can be expressed in terms of 
Hirzebruch L-functions \cite{zagier}. 
\end{theorem} 
It should be mentioned that a version of the higher-dimensional Dedekind sums had already been introduced by Carlitz 
\cite{carlitz}: 
  \[ \sum_{ k_{1}, \dots, k_{d} \text{ mod } a_{0} } \left( \left( \frac{ a_{1} k_{1} + \dots + a_{d} k_{d} }{ a_{0} } \right) \right) \left( \left( \frac{ a_{1} }{ a_{0} } \right) \right) \cdots \left( \left( \frac{ a_{d} }{ a_{0} } \right) \right) \ . \] 
Berndt \cite{berndt} noticed that these sums are, up to trivial factor, Zagier's higher-dimensional Dedekind sums. 

If we write the higher-dimensional Dedekind sum as a sum over roots of unity, 
  \[ \s ( a_{0} ; a_{1} , \dots , a_{d} ) = \frac{ 1 }{ a_{0} } \sum_{ \lambda^{ a_{0} } = 1 \not= \lambda } \frac{ \lambda^{ a_{1} } + 1 }{ \lambda^{ a_{1} } - 1 } \cdots \frac{ \lambda^{ a_{d} } + 1 }{ \lambda^{ a_{d} } - 1 } \ , \] 
it becomes clear that it suffices to study sums of the form 
  \[ \frac{ 1 }{ a_{0} } \sum_{ \lambda^{ a_{0} } = 1 \not= \lambda } \frac{ 1 }{ ( \lambda^{ a_{1} } - 1 ) \cdots ( \lambda^{ a_{d} } - 1 ) } \ . \] 
Zagier's Dedekind sum can be expressed as a sum of expressions of this kind. 
On the other hand, we consider special cases of the Dedekind-Rademacher sum, namely, for $ n \in \Z $, 
  \[ \s \left(a,b; \frac{ n }{ b } , 0 \right) = \sum_{ k \text{ mod } b } \left( \left( \frac{ ka + n }{ b } \right) \right) \left( \left( \frac{ k }{ b } \right) \right) \ . \] 
Knuth \cite{knuth} discovered that these generalized Dedekind sums describe 
the statistics of pseudo random number generators. 
In \cite{br}, we used the convolution theorem for finite Fourier series to 
show that, if $a$ and $b$ are relatively prime, 
  \begin{equation}\label{raddedsum} \s \left(a,b; \frac{ n }{ b } , 0 \right) = - \frac{ 1 }{ b } \sum_{ \lambda^{ b } = 1 \not= \lambda } \frac{ \lambda^{ -n } }{ ( 1 - \lambda^{ a } ) ( 1 - \lambda ) } - \frac{1}{2} \left\{ \frac{ n }{ b } \right\} + \frac{1}{4} - \frac{ 1 }{ 4b } \ . \end{equation} 
Here $ \{ x \} = x - [x] $ denotes the fractional part of $x$. 
Comparing this with the representation we obtained for Zagier's Dedekind sums motivates the study of the Fourier-Dedekind sum 
  \[ \sigma_{n} \left( a_{ 1 } , \dots, a_{ d } ; a_{0} \right) = \frac{1}{a_{0}} \sum_{ \lambda^{ a_{0} } = 1 } \frac{ \lambda^{ n } }{ \left( 1 - \lambda^{ a_{ 1 } } \right) \cdots \left( 1 - \lambda^{ a_{ d } } \right) } \ , \]
a finite Fourier series in $n$. 
Gessel \cite{gessel} gave a new reciprocity law for a special case of Fourier-Dedekind sums:
\begin{theorem}[Gessel]\label{ges} Let $p$ and $q$ be relatively prime and suppose that $ 1 \leq n \leq p+q $. Then 
  \begin{eqnarray*} &\mbox{}& \frac{ 1 }{ p }  \sum_{ \lambda^{ p } = 1 \not= \lambda } \frac{ \lambda^{ n }  }{ \left( 1 - \lambda^{ q } \right) \left( 1 - \lambda \right)  } +  \frac{ 1 }{ q }  \sum_{ \lambda^{ q } = 1 \not= \lambda } \frac{ \lambda^{ n }  }{ \left( 1 - \lambda^{ p  }  \right) \left( 1 - \lambda \right)  } \\ 
                    &\mbox{}& \qquad = - \frac{ n^{2} }{ 2pq } + \frac{n}{2} \left( \frac{1}{p} + \frac{1}{q} + \frac{1}{pq} \right) - \frac{1}{4} \left( \frac{1}{p} + \frac{1}{q} + 1 \right)  - \frac{1}{12} \left( \frac{p}{q} + \frac{1}{pq} + \frac{q}{p} \right) \ . \end{eqnarray*} 
\end{theorem} 
It is easy to see that the reciprocity law for classical Dedekind sums (Theorem \ref{dedrec}) is a special 
case of Gessel's theorem. We can rephrase the statement of Gessel's theorem in terms of Dedekind-Rademacher sums by means of (\ref{raddedsum}): 
for $p$ and $q$ relatively prime, and $ 1 \leq n \leq p+q $, 
  \begin{align*} &\s \left( q, p; \frac{ -n }{ p } , 0 \right) + \s \left( p, q; \frac{ -n }{ q } , 0 \right) \\ 
                 &\stackrel{ \text{def} }{ = } \ \sum_{ k=0 }^{ p-1 } \left( \left( \frac{ q k - n }{ p } \right) \right) \left( \left( \frac{ k }{ p } \right) \right) + \sum_{ k=0 }^{ q-1 } \left( \left( \frac{ p k - n }{ q } \right) \right) \left( \left( \frac{ k }{ q } \right) \right) \\ 
                 &= \frac{ n^{2} }{ 2pq } - \frac{n}{2} \left( \frac{1}{p} + \frac{1}{q} + \frac{1}{pq} \right) + \frac{1}{4} + \frac{1}{12} \left( \frac{p}{q} + \frac{1}{pq} + \frac{q}{p} \right) - \frac{1}{2} \left\{ \frac{ -t }{ p } \right\} - \frac{1}{2} \left\{ \frac{ -t }{ q } \right\} \ . \end{align*} 
We will now view the Fourier-Dedekind sum from a generating-function point of view, 
which will allow us to obtain and extend geometric proofs of Dedekind's, 
Zagier's and Gessel's reciprocity laws. 


\section{A new combinatorial identity for partitions from a finite set}\label{comb} 
The form of the Fourier-Dedekind sum 
  \[ \sigma_{-n} \left( a_{ 1 } , \dots, a_{ d } ; a_{ 0 } \right) = \frac{1}{ a_{ 0 } } \sum_{ \lambda^{ a_{ 0 } } = 1 \not= \lambda } \frac{ \lambda^{ -n } }{ \left( 1 - \lambda^{ a_{ 1 } } \right) \cdots \left( 1 - \lambda^{ a_{ d } } \right) } \]
suggests the use of a generating function 
  \[ f(z) := \frac{ 1 }{ 1 - z^{ a_{ 0 } } } \ \frac{ z^{ -n } }{ \left( 1 - z^{ a_{ 1 } } \right) \cdots \left( 1 - z^{ a_{ d } } \right) } \ . \]
In fact, let's expand this generating function into partial fractions: suppose, for simplicity, that $ n > 0 $, and 
$ a_{ 0 } , \dots, a_{ d } $ are pairwise relatively prime. Then we can write 
  \[ f(z) = \sum_{ \lambda^{ a_{ 0 } } = 1 \not= \lambda } \frac{ A_{ \lambda } }{ z - \lambda } + \dots + \sum_{ \lambda^{ a_{ d } } = 1 \not= \lambda } \frac{ A_{ \lambda } }{ z - \lambda } + \sum_{ k=1 }^{ d+1 } \frac{ B_{k} }{ (z-1)^{ k } } + \sum_{ k=1 }^{ n } \frac{ C_{k} }{ z^{ k } } \ . \] 
The coefficient $ A_{ \lambda } $ for, say, a nontrivial $ a_{0} $'th root of unity $ \lambda $ can be derived easily: 
  \[ A_{ \lambda } = \lim_{ z \to \lambda } ( z - \lambda ) f(z) = - \frac{ \lambda }{ a_{ 0 } } \frac{ \lambda^{ -n } }{ \left( 1 - \lambda^{ a_{ 1 } } \right) \cdots \left( 1 - \lambda^{ a_{ n } } \right) } \ . \] 
Hence we obtain the Fourier-Dedekind sums if we consider the constant coefficient of $f$ (in the Laurent series about $ z=0 $):
  \begin{eqnarray} &\mbox{}& \const (f) = \sum_{ \lambda^{ a_{ 0 } } = 1 \not= \lambda } \frac{ A_{ \lambda } }{ - \lambda } + \dots + \sum_{ \lambda^{ a_{ d } } = 1 \not= \lambda } \frac{ A_{ \lambda } }{ - \lambda } + \sum_{ k=1 }^{ d+1 } (-1)^{ k } B_{k} \label{part} \\ 
                   &\mbox{}& \qquad = \sigma_{-n} \left( a_{ 1 } , \dots, a_{ d } ; a_{ 0 } \right) + \dots + \sigma_{-n} \left( a_{ 0 } , \dots, a_{ d-1 } ; a_{ n } \right) + \sum_{ k=1 }^{ d+1 } (-1)^{ k } B_{k} \ . \nonumber \end{eqnarray} 
The coefficients $ B_{k} $ are simply the coefficients of the Laurent series of $f$ about $ z=1 $, and are easily
computed, by hand or using mathematics software such as {\tt Maple} or {\tt Mathematica}. It is not hard to see 
that they are polynomials in $n$ whose coefficients are rational functions of the 
$ a_{ 0 } , \dots, a_{ d } $.\footnote{After this paper was submitted, general formulas for these polynomials were discovered in \cite{bgk}.}\label{fn} 
To simplify notation, define
  \begin{equation}\label{pol} q ( a_{ 0 } , \dots, a_{ d }, n ) := \sum_{ k=1 }^{ d+1 } (-1)^{ k } B_{k} \ . \end{equation}

On the other hand, we can compute the constant coefficient of $f$ by brute force: By expanding 
  \[ f(z) = \left( \sum_{ k_{0} \geq 0 } z^{ k_{0} a_{0} } \right) \cdots \left( \sum_{ k_{d} \geq 0 } z^{ k_{d} a_{d} } \right) z^{ -n } \ , \] 
we can see that $\const(f)$ enumerates the ways of writing $n$ as a linear combination of the $ a_{ 0 } , \dots, a_{ d } $ with nonnegative coefficients: 
  \begin{align} \const (f) &= \# \left\{ \left( k_{0} , \dots, k_{d} \right) \in \Z^{ d+1 } : \ k_{j} \geq 0 , \ k_{0} a_{0} + \dots + k_{d} a_{d} = n \right\} \label{triv} \\ 
                           &= p_{ \{ a_{ 0 } , \dots, a_{ d } \} } (n) \ . \nonumber \end{align} 
This defines the partition function with parts in the finite set 
$ A := \{ a_{ 0 } , \dots, a_{ n } \} $. Geometrically, 
$ p_{A} (n) $ enumerates the integer points in $n$-dilates of the rational polytope
  \[ \P := \left\{ \left( x_{0} , \dots, x_{d} \right) \in \R^{ d+1 } : \ x_{j} \geq 0 , \ x_{0} a_{0} + \dots + x_{d} a_{d} = 1 \right\} \ . \]
This geometric interpretation allows us to use the machinery of Ehrhart theory, which 
will be advantageous in the following section.  We next give an explicit formula
for the famous `coin-exchange problem'---that is, the number of ways to form $n$ cents
from a finite set of coins with given denominations $ a_{ 0 } , \dots, a_{ d } $: comparing (\ref{part}) with (\ref{triv}) yields our central result \cite{bdr}. 
\begin{theorem}\label{main}
Suppose $ a_{ 0 } , \dots, a_{ d } $ are pairwise 
relatively prime, positive integers.  We recall that the number of partitions of an
integer $n$ from the finite set of $a_i$'s is defined by
  \[ p_{ \{ a_{ 0 } , \dots, a_{ d } \} } (n) := \left\{ \left( k_{0} , \dots, k_{d} \right) \in \Z^{ d+1 } : \ k_{j} \geq 0 , \ k_{0} a_{0} + \dots + k_{d} a_{d} = n \right\} . \] 
Then   
  \[ p_{ \{ a_{ 0 } , \dots, a_{ d } \} } (n) = q ( a_{ 0 } , \dots, a_{ d }, n ) + \sum_{ j=0 }^{ d } \sigma_{-n} \left( a_{ 0 } , \dots, \hat{ a_{ j } }, \dots, a_{ d } ; a_{ j } \right) \ , \]
where $ q ( a_{ 0 } , \dots, a_{ d }, n ) $ is given by {\rm (\ref{pol})}.
\end{theorem} 
The first few expressions for $ q ( a_{ 0 } , \dots, a_{ d }, n ) $ are 
\begin{eqnarray*} &\mbox{}& q ( a_{ 0 }, n ) = \frac{ 1 }{ a_{ 0 } } \\ 
                  &\mbox{}& q ( a_{ 0 }, a_{ 1 }, n ) = \frac{ n }{ a_{ 0 } a_{ 1 } } + \frac{ 1 }{ 2 } \left( \frac{ 1 }{ a_{ 0 } } + \frac{ 1 }{ a_{ 1 } }  \right) \\
                  &\mbox{}& q ( a_{ 0 } , a_{ 1 }, a_{ 2 }, n ) = \frac{ n^{ 2 }  }{ 2 a_{ 0 } a_{ 1 } a_{ 2 }  } + \frac{ n }{ 2 } \left( \frac{ 1 }{ a_{ 0 } a_{ 1 }  } + \frac{ 1 }{ a_{ 0 } a_{ 2 }  } + \frac{ 1 }{ a_{ 1 } a_{ 2 }  }  \right) \\
                  &\mbox{}& \qquad + \frac{ 1 }{ 12 } \left( \frac{ 3 }{ a_{ 0 }  } + \frac{ 3 }{ a_{ 1 }  } + \frac{ 3 }{ a_{ 2 }  } + \frac{ a_{ 0 }  }{ a_{ 1 } a_{ 2 }  } + \frac{ a_{ 1 }  }{ a_{ 0 } a_{ 2 }  } + \frac{ a_{ 2 }  }{ a_{ 0 } a_{ 1 }  }  \right) \\
                  &\mbox{}& q ( a_{ 0 } , a_{ 1 }, a_{ 2 }, a_{ 3 }, n ) = \frac{ n^{ 3 }  }{ 6 a_{ 0 } a_{ 1 } a_{ 2 } a_{ 3 }  } \\ 
                  &\mbox{}& \qquad + \frac{ n^{ 2 }  }{ 4 } \left( \frac{ 1 }{ a_{ 0 } a_{ 1 } a_{ 2 }  } + \frac{ 1 }{ a_{ 0 } a_{ 1 } a_{ 3 }  } + \frac{ 1 }{ a_{ 0 } a_{ 2 } a_{ 3 }  } + \frac{ 1 }{ a_{ 1 } a_{ 2 } a_{ 3 }  }  \right) \\
                  &\mbox{}& \qquad + \frac{ n }{ 4 } \left( \frac{ 1 }{ a_{ 0 } a_{ 1 }  } + \frac{ 1 }{ a_{ 0 } a_{ 2 }  } + \frac{ 1 }{ a_{ 0 } a_{ 3 } } + \frac{ 1 }{ a_{ 1 } a_{ 2 }  } + \frac{ 1 }{ a_{ 1 } a_{ 3 }  } + \frac{ 1 }{ a_{ 2 } a_{ 3 }  } \right) \\
                  &\mbox{}& \qquad + \frac{ n }{ 12 } \left( \frac{ a_{ 0 }  }{ a_{ 1 } a_{ 2 } a_{ 3 }  } + \frac{ a_{ 1 }  }{ a_{ 0 } a_{ 2 } a_{ 3 }  } + \frac{ a_{ 2 }  }{ a_{ 0 } a_{ 1 } a_{ 3 }  } + \frac{ a_{ 3 }  }{ a_{ 0 } a_{ 1 } a_{ 2 }  }  \right) \\
                  &\mbox{}& \qquad + \frac{ 1 }{ 24 } \left( \frac{ a_{ 0 }  }{ a_{ 1 } a_{ 2 }  } + \frac{ a_{ 0 }  }{ a_{ 1 } a_{ 3 }  } + \frac{ a_{ 0 }  }{ a_{ 2 } a_{ 3 }  } + \frac{ a_{ 1 }  }{ a_{ 0 } a_{ 2 }  } + \frac{ a_{ 1 }  }{ a_{ 0 } a_{ 3 }  } + \frac{ a_{ 1 }  }{ a_{ 2 } a_{ 3 }  } \right. \\
                  &\mbox{}& \qquad \quad \left. + \frac{ a_{ 2 }  }{ a_{ 0 } a_{ 1 }  } + \frac{ a_{ 2 }  }{ a_{ 0 } a_{ 3 }  } + \frac{ a_{ 2 }  }{ a_{ 1 } a_{ 3 }  } + \frac{ a_{ 3 }  }{ a_{ 0 } a_{ 1 }  } + \frac{ a_{ 3 }  }{ a_{ 0 } a_{ 2 }  } + \frac{ a_{ 3 }  }{ a_{ 1 } a_{ 2 }  } \right) \\
                  &\mbox{}& \qquad + \frac{ 1 }{ 8 } \left( \frac{ 1 }{ a_{ 0 }  } + \frac{ 1 }{ a_{ 1 }  } + \frac{ 1 }{ a_{ 2 }  } + \frac{ 1 }{ a_{ 3 }  } \right) \ . \end{eqnarray*} 


\section{Reciprocity laws}\label{recsec} 
We will now use Theorem \ref{main} to prove and extend some of the reciprocity theorems stated earlier.
We will make use of two results due to Ehrhart for {\bf rational polytopes}, that is, polytopes whose
vertices are rational. Ehrhart \cite{ehrhart} initiated the study of the number of
integer points (``lattice points'') in integer dilates of such polytopes:
\begin{definition} Let $ \P \subset \R^{d} $ be a rational polytope, and $n$ a positive integer. We denote the
number of lattice points in the dilates of the closure of $\P$ and its interior by
  \[ L ( \overline{\P} , n ) := \# \left( n \P \cap \Z^{d} \right) \qquad \mbox{ and } \qquad L ( {\P}^{\circ} , n ) := \# \left( n \P^{\circ} \cap \Z^{d} \right) \ , \]
respectively.
\end{definition}
Ehrhart proved that $ L ( \overline{\P} , n ) $ and $ L ( {\P}^{\circ} , n ) $ are
{\bf quasipolynomials} in the integer variable $n$, that is, expressions of the form
  \[ c_{d}(n) \ n^{d} + \dots + c_{1}(n) \ n + c_{0}(n) \ , \]
where $c_{0}, \dots , c_{d}$ are periodic functions in $n$.
Ehrhart conjectured the following fundamental theorem, which establishes an algebraic connection
between our two lattice-point-count operators. Its original proof is due to Macdonald \cite{macdonald}. 
\begin{theorem}[Ehrhart-Macdonald reciprocity law]\label{ehrrec} Suppose the rational po\-lytope $ \P $ is
homeomorphic to a $d$-manifold. Then
  \[ L ( {\P}^{\circ}  , -n ) = (-1)^{ d } L ( \overline{\P} , n )  \ . \]
\end{theorem}
This enables us to rephrase Theorem \ref{main} for the quantity
  \[ p_{ \{ a_{ 0 } , \dots, a_{ d } \} }^{ \circ } (n) := \# \left\{ \left( k_{0} , \dots, k_{d} \right) \in \Z^{ d+1 } : \ k_{j} > 0 , \ k_{0} a_{0} + \dots + k_{d} a_{d} = n \right\} \ . \]
By Theorem \ref{ehrrec}, we have the following result. 
\begin{corollary}\label{maincor} Suppose $ a_{ 0 } , \dots, a_{ d } \in \N $ are pairwise relatively prime. Then
  \[ p_{ \{ a_{ 0 } , \dots, a_{ d } \} }^{ \circ } (n) = (-1)^{ d } \left( q ( a_{ 0 } , \dots, a_{ d }, -n ) + \sum_{ j=0 }^{ d } \sigma_{n} \left( a_{ 0 } , \dots, \hat{ a_{ j } }, \dots, a_{ d } ; a_{ j } \right) \right) \ , \]
where $ q ( a_{ 0 } , \dots, a_{ d }, n ) $ is given by {\rm (\ref{pol})}.
\end{corollary}
We note that we could have derived this identity from scratch in a similar way as Theorem \ref{main}, without using
Ehrhart-Macdonald reciprocity.

The reason for switching to $ p_{ \{ a_{ 0 } , \dots, a_{ d } \} }^{ \circ } (n) $ is that
  \[ p_{ \{ a_{ 0 } , \dots, a_{ d } \} }^{ \circ } (n) = 0 \]
for $ 0 < n < a_{ 0 } + \dots + a_{ d } $, by the very definition of $ p_{ \{ a_{ 0 } , \dots, a_{ d } \} }^{ \circ } (n) $.
This yields a reciprocity law:
\begin{theorem} Let $ a_{0} , \dots , a_{d} $ be pairwise relatively prime integers and $0 < n < a_{ 0 } + \dots + a_{ d }$. Then
  \[ \sum_{ j=0 }^{ d } \sigma_{n} ( a_{ 0 } , \dots, {\hat a_{ j } } , \dots , a_{ d } ; a_{ j } )  = - q ( a_{ 0 } , \dots, a_{ d }, -n )  \ , \]
where $ q ( a_{ 0 } , \dots, a_{ d }, n ) $ is given by {\rm (\ref{pol})}.
\end{theorem}
For $ d = 2, a_{0} = p, a_{1} = q, a_{2} = 1 $, this is the statement of Gessel's Theorem \ref{ges}, which, in turn, implies
Dedekind's reciprocity law Theorem \ref{dedrec}.

To prove Zagier's Theorem \ref{zagierthm} in the language of Fourier-Dedekind sums, we make use of another result of Ehrhart \cite{ehrhart} 
on {\bf lattice polytopes}, that is, polytopes whose vertices have integer coordinates.
Recall that the {\bf reduced Euler characteristic} of a polytope $\P$ can be defined as
  \[ \chi ( \P ) := \sum_{ \sigma } (-1)^{ \dim \sigma } \ , \] 
where the sum is over all sub-simplices of $\P$.
\begin{theorem}[Ehrhart]\label{ehrthm} Let $\P$ be a lattice polytope. Then $ L ( \overline{\P} , n ) $
is a polynomial in $n$ whose constant term is $ \chi ( \P ) $.
\end{theorem}
We note that the polytope $\P$ corresponding to $ p_{ \{ a_{ 0 } , \dots, a_{ d } \} } (n) $ is convex and
hence has Euler characteristic 1. If we now dilate $\P$ only by multiples of $ a_{ 0 } \cdots a_{ d } $,
say $ n = a_{ 0 } \cdots a_{ d } m $, we obtain the dilates of a lattice polytope. Theorem \ref{main}
simplifies for these $n$ to
  \[ p_{ \{ a_{ 0 } , \dots, a_{ d } \} } ( a_{ 0 } \cdots a_{ d } m ) = q ( a_{ 0 } , \dots, a_{ d }, a_{ 0 } \cdots a_{ d } m ) + \sum_{ j=0 }^{ d } \sigma_{0} \left( a_{ 0 } , \dots, \hat{ a_{ j } }, \dots, a_{ d } ; a_{ j } \right) \ , \]
by the periodicity of the Fourier-Dedekind sums. However, $\chi(\P)=1$, and Theorem \ref{ehrthm} yields
a result equivalent to Zagier's reciprocity law for his higher-dimensional Dedekind sums, Theorem \ref{zagierthm}:
\begin{theorem}\label{zag} For pairwise relatively prime integers $ a_{1} , \dots , a_{d} $,
  \[ \sum_{ j=0 }^{ d } \sigma_{0} ( a_{ 0 } , \dots, {\hat a_{ j } } , \dots , a_{ d } ; a_{ j } ) = 1 - q ( a_{ 0 } , \dots, a_{ d }, 0 ) \ , \]
where $ q ( a_{ 0 } , \dots, a_{ d }, n ) $ is given by {\rm (\ref{pol})}.
\end{theorem}


\section{The computational complexity of Zagier's \\ higher-dimensional Dedekind sums}\label{last}

In this section we give a proof of the polynomial-time complexity of Zagier's
higher-dimensional Dedekind sums, in fixed dimension $d$.
In \cite{barv}, there is a nice theorem due to Barvinok which guarantees 
the polynomial-time
computability of the generating function attached to a rational polyhedron.
We will use his  theorem for a cone.  First, we mention that a common way to 
enumerate lattice points in a cone $\K$ (and in polytopes) is to use 
the generating function
  \[ f(\K, \x) := \sum_{\m \in \K \cap \Z^d} \x^{\m} \ , \] 
where we use the standard multivariate notation $\x^{\m} := x_1^{m_1} \dots x_d^{m_d}$.
It is an elementary fact that for rational cones these generating functions
are always rational functions of the variable $\x$. Barvinok's theorem reads as follows:

\begin{theorem}[Barvinok] Let us fix the dimension $d$. There exists a polyno\-mial-time algorithm, which 
for a given rational polyhedron $\K \subset \R^d$, 
  \[ \K = \{ \x \in \R^{d} : \langle \c_{i}, \x \rangle \leq \beta_{i}, i = 1 ... m \} , \ where   \ \c_i \in \Z^d  \ and  \ \beta_i \in \Q \] 
computes the generating function
  \[ f(\K, \x) := \sum_{\m \in \K \cap \Z^d} \x^{\m} \] 
in the form (a virtual decomposition) 
  \[ f(\K, \x) = \sum_{i \in I}\epsilon_i \frac{{\x}^{a_i}}{  (1-\x^{b_{i1}}) \cdots (1-\x^{b_{id}})} \ , \] 
where $\epsilon_i \in \{ -1, 1 \}, a_i \in \Z$, and $b_{i1}, ... , b_{id}$ is a basis of $\Z^d$ for each $i$. 
The computational complexity of the algorithm for finding this virtual decomposition is $\mathcal L^{O(d)}$, where $\mathcal L$ 
is the input size of $\K$. In particular, the number $I$ of terms in the summand is $\mathcal L^{O(d)}$. 
\end{theorem}

Thus Barvinok's algorithm finds the coefficients of the rational function $f(\K, \x)$ in polynomial time. 
In \cite{diaz}, on the other hand, the generating function $f(\K, \x)$ is given in terms of an average over a finite abelian group 
of  a product of $d$  cotangent functions, whose arguments are in terms of the coordinates of the vertices which generate the cone $\K$ (these are 
the extreme points of $\K$ whose convex hull is $\K$).  This is the main theorem in 
\cite{diaz} and we apply it below to a special lattice cone which will give us the Zagier-Dedekind sums we want to study. 

The following theorem is part of a bigger project on the computability 
of generalized Dedekind sums in all dimensions. A slightly different proof is sketched in \cite{barv}.
\begin{theorem} For fixed dimension $d$, the higher-dimensional Dedekind sums 
  \[ \s ( a_{0} ; a_{1} , \dots , a_{d} ) = \frac{ (-1)^{d/2} }{ a_{0} } \sum_{ k=1 }^{ a_{0} - 1 } \cot \frac{ \pi k a_{1} }{ a_{0} } \cdots \cot \frac{ \pi k a_{d} }{ a_{0} } \ . \]
are polynomial-time computable. 
\end{theorem} 
\begin{proof} 
Let $\K \subset \R^{d+1}$ be the cone generated by the positive real span of the vectors 
\begin{eqnarray*} v_1 &=& (1, 0, \dots , 0, a_1) \\
                  v_2 &=& (0, 1, 0, \dots , 0 , a_2) \\
                  & \vdots \\
                  v_{d} &=& (0, \dots , 0, 1, a_{d}) \\
                  v_{d+1} &=& (0, \dots , 0, a_0) \end{eqnarray*}
Then the right-hand side of the main theorem of \cite{diaz}
is in this case
$$
\frac{1}{2^{d+1} a_0}
\sum_{ k=1 }^{ a_{0} - 1 } \prod_{j=0}^{d} \left( 1 + \coth \frac{\pi}{a_{0}}\left( s + i a_{j} k \right) \right).
$$
When we compute the coefficient of $s^{-1}$ in this meromorphic function of $s$, we arrive at the
following higher-dimensional Dedekind sum:
  \[ \sum_{ k=1 }^{ a_{0} - 1 } \cot \frac{ \pi k a_{1} }{ a_{0} } \cdots \cot \frac{ \pi k a_{d} }{ a_{0} } \]
plus other products of lower-dimensional Zagier-Dedekind sums.
By induction on the dimension, all of the lower-dimensional Zagier-Dedekind
sums are polynomial-time computable, and since the left-hand
side of the main theorem is polynomial-time computable by Barvinok's
theorem, the above Zagier-Dedekind sums in dimension $d$ is now also
polynomial-time computable.  
\end{proof}


\bibliographystyle{amsalpha}

\end{document}